\documentclass[12pt]{amsart}

\usepackage{txfonts}
\usepackage{amscd}
\usepackage{cite}
\usepackage{mathbbold, color}
\usepackage{mathbbold, cite}
\usepackage{epsfig}
\usepackage{verbatim}
\usepackage{mathdots}
\usepackage{amssymb}
\usepackage{amsfonts}
\usepackage{amsbsy}
\usepackage{graphicx}
\usepackage{enumerate}
\usepackage{ulem}

\newtheorem{theo}{Theorem}[section]
\newtheorem{lem}[theo]{Lemma}
\newtheorem{prop}[theo]{Proposition}

\theoremstyle{definition}
\newtheorem{defi}[theo]{Definition}

\newtheorem{rem}[theo]{Remark}

\numberwithin{equation}{section}

\newcommand{\R}{{\mathbb R}}
\newcommand{\Z}{{\mathbb Z}}

\newcommand{\N}{{\mathbb N}}

\newcommand{\D}{{\mathcal D}}

\newcommand{\eproof}{\hfill$\square$}

\begin{document}

\nocite{*}

        \baselineskip 16pt
        
        \title[]{Beurling densities of regular maximal orthogonal sets of self-similar spectral measure with consecutive digit sets}
        
        \author[Yu-Liang Wu]{Yu-Liang Wu}
        \address[Yu-Liang Wu]{Department of Mathematical Sciences\\
                University of Oulu, P.O. Box 3000, 90014, Oulu, Finland}
        \email{yu-liang.wu@oulu.fi}
        \author[Zhi-Yi Wu]{Zhi-Yi Wu}
        \address[Zhi-Yi Wu]{School of Mathematics and Information Science, Guangzhou University, Guang\-zhou, 510006, P.~R.~China;~Department of Mathematical Sciences\\
                University of Oulu, P.O. Box 3000, 90014, Oulu, Finland}
        \email{zhiyiwu@126.com}
        
%
%
        %
        
        
        
\begin{abstract}
Beurling density plays a key role in the study of frame-spectrality of normalized Lebesgue measure restricted to a set. Accordingly, in this paper, the authors study the $s$-Beurling densities of regular maximal orthogonal sets of a class of self-similar spectral measures, where $s$ is the Hausdorff dimension of its support and obtain their exact upper bound of the
densities.

\end{abstract}
        
        \maketitle
        \section{Introduction}
\begin{defi}
Let $\mu$ be a Borel probability measure with compact support in $\R^d.$ We say that $\{e^{-2\pi i\lambda x}\}_{\lambda\in\Lambda}$  is a \textit{Fourier frame} of the Hilbert space $L^2(\mu)$ if there exist two constants $A, B>0$ and  a countable set $\Lambda\subseteq \R^d$, called a \textit{frame spectrum},  such that for every $f\in L^2(\mu)$, we have
\begin{equation}\label{frame}
A\|f\|^2\le \sum_{\lambda\in \Lambda}|\langle f, e_\lambda \rangle|^2\le B\|f\|^2,
\end{equation}
where $\langle \cdot, \cdot\rangle $ is the inner product in $L^2(\mu)$, $e_\lambda(x)=e^{2\pi i x\cdot\lambda}$ and $x\cdot\lambda=\sum_{i=1}^dx_i\lambda_i$ is the standard inner product in $\R^d$. In this case, we call $\mu$  a {\it frame spectral measure}. Specially, if $A=B=1$ in equation $\eqref{frame}$, it is easy to see that the set $\{e_{\lambda}\}_{\lambda\in\Lambda}$ is an orthonormal basis for $L^2(\mu)$. In this case, we call $\mu$ a {\it spectral measure} and the set $\Lambda$ a {\it spectrum}.
\end{defi}

It is known that a frame-spectral measure is either a finite discrete measure, absolutely continuous or singular continuous measure with respect to Lebesgue one \cite{HLL13}. When $\mu$ is the Lebesgue measure supported on $[0,1]$, the work of Landau, Jaffard, and Seip relates the frame spectrum of $L|_{[0,1]}$ closely with the Beurling densities.

\begin{defi}Let $\Lambda$ be a countable set in $\R$. The {\it upper and lower Beurling density} of $\Lambda$ are defined, respectively, by
$$ D^{+}(\Lambda)=\limsup\limits_{h\rightarrow\infty}\sup\limits_{x\in\R}\frac{\#(\Lambda\cap (x-h,x+h))}{h}$$
 and
$$ D^{-}(\Lambda)=\liminf\limits_{h\rightarrow\infty}\inf\limits_{x\in\R}\frac{\#(\Lambda\cap (x-h,x+h))}{h},$$
where $\#E$ is the cardinality of the set $E$.
\end{defi}
Using the notion of Beurling density, the frame properties of $\{e^{2\pi i\lambda x}\}_{\lambda\in\Lambda}$ are almost characterized \cite{Lan67,Ja91,Sei95}:
\begin{theo}
For $\{e^{-2\pi i\lambda x}\}_{\lambda\in\Lambda}$ to be a Fourier frame for $L^2([0,1])$, it is necessary that $D^{+}(\Lambda)<\infty$ and $D^{-}(\Lambda)\geq 2$, and it is sufficient that $D^{+}(\Lambda)<\infty$ and $D^{-}(\Lambda)>2$.
\end{theo}

A complete description of the problem that those sequences $\Lambda$ generate Fourier frames is obtained by Ortega-Cerd\`{a} and Seip \cite{OS02}, who solved the critical case when $D^{-}(\Lambda)=2$ by using de Branges' theory of Hilbert space of entire functions.

For absolutely continuous measures $\mathrm{d}\mu=g(x)\mathrm{d}x$, Lai \cite{Lai11} proved that  if there exists a Fourier frame, then the function $g$ must be bounded above and below on its support. The Beurling density also plays a key role in his proof. In this paper, we focus on the singular case, which is muss less understood.

In 1998, Jorgensen and Pederson \cite{JP} discovered that the standard middle-fourth Cantor measure $\mu_{4,\{0,2\}}$ is a spectral measure, which is the first non-atomic and singular spectral measure. In the same paper, they proved the standard middle-third Cantor measure $\mu_{3,\{0,2\}}$ is not spectral. Since then, the spectrality of self-similar measures/self-affine measures/Moran measures and related properties have attracted a
        great amount of attention (see \cite{LWY,AW,D,D16,DLau,DH,DHS,LJJ15,DHL13,AH,AFL,DFY,DJ07,DL14,DHLai19,HTW19} and the references therein for
        more details). Meanwhile, many new phenomena not previously observed on the Lebesgue measure are unveiled\cite{LWY,FHW18,HKTW,HTW19,DHS14,Str00,Str06,TW}. In particular, {\L}aba and Wang \cite{LWY} discovered that a singularly continuous self-similar spectral measure may admit many spectra. The exotic phenomenon naturally  lead researchers to do more sophisticated analysis of the structure of the spectra \cite{D16,DHL13,DHS,FHW18,HTW19,LJJ15}.


To completely characterize the spectra for a self-similar spectral measure, Dutkay et al. \cite{DHSW11} use the {\it Beurling dimension} as a replacement for Beurling density. Let $\Lambda\subseteq \R^d$ be a countable set. For $r>0,$ the \textit{upper $r$-Beurling density} of $\Lambda$ is defined by
\[
D_r^+(\Lambda)=\limsup\limits_{h\to \infty}\sup_{x\in \R^d}\frac{\#(\Lambda \cap B(x,h))}{h^r},
\]
where  $B(x,h)$ is the open ball centered at $x$ with radius $h$.
It is not difficult to prove that there is a critical value of $r$ where $D_r^+(\Lambda)$ jumps from $\infty$ to $0$. This critical value is defined as the \textit{Beurling dimension} of $\Lambda$, or more precisely,
\[
\dim_{Be} (\Lambda)=\inf\{r:D_r^+(\Lambda)=0\}=\sup\{r:D_r^+(\Lambda)=\infty\}.
\]
Duktay et al. \cite {DHSW11} proved that the Beurling dimension of the spectra of self-similar spectral measure which satisfies the open set condition is not greater than the Hausdorff dimension of its support, and they may equal under some mild condition. On the other hand, there exist arbitrarily sparse spectra for many singular spectral measures \cite{AL,DHL13}, i.e., there exist spectra with Beurling dimensions zero, which is in stark contrast with the case of the Lebesgue measure.

However, for a long time, in contrast to the classic case, the notion of Beurling density did not seem to play an important role in fractal spectral measure theory. Motivated by this, in this paper, we investigate the Beurling density of the maximal orthogonal sets of a class of self-similar spectral measure. Our ultimate goal is that we want to know how the values of Beurling densities of maximal orthogonal sets imply the completeness.

Now, we recall the definition of self-similar measures. A {\it self-similar} measure on $\R^d$ is defined by Hutchinson \cite{Hut} to be a probability measure $\mu:=\mu_{A,\D}$ satisfying
        \begin{equation}\label{eqselfmea}
                \mu=\frac{1}{\#\D}\sum_{d\in\D}\mu\circ f_d^{-1},
        \end{equation}
        where $\D\subseteq\R^d$ is a finite digit set and $\{f_d=A^{-1}(x +d)\}_{d\in\D}$ is a sequence of contractive similarity with $A=rQ$ ($r>1$ and $Q$ orthogonal). Here, $\#A$ is the cardinality of $A$. A noteworthy expression for such measures is given by the infinite convolution product
        \begin{equation}\label{eqconv}
                \mu_{A,\D}=\delta_{A^{-1}\D}\ast\delta_{A^{-2}\D}\ast\cdots,
        \end{equation}
        where $\delta_E=\frac{1}{\#E}\sum_{e\in E}\delta_e$, $\delta_e$ is the Dirac measure at the point $e\in E$ and the convergence is in weak sense.

The following discussions will focus on self-similar measures on $\R$, where $A=p$ be a real number strictly  larger than $1$ and $\D=r\{0,1,\ldots,q-1\}$ with $r=\frac{p}{q}$. In this case, we call $\mu_{p,q}:=\mu_{p,\D}$ a self-similar measure with consecutive digit set. When $q=2$, it reduces to the Bernoulli convolutions. The spectrality of the  Bernoulli convolution was considered by Hu and Lau \cite{HL}, and was completely resolved by Dai \cite{D}. They proved that $\mu_{p,2}$ is a spectral measure if and only if $p \in 2\Z$. These results are generalized further to the $\mu_{p,q}$ with $p>2$ by Dai, He and Lai \cite{DHL14}, in which they proved that $\mu_{p,q}$ is a spectral measure if and only if $p\in\Z$ and $q$ divides $p$. In this case, the  structure of maximal orthogonal sets for $\mu_{p,q}$ is well characterized \cite{DHL13}, and our results are based on their constructions, as described below.
        
        
        For $q\ge 2,$ write $\Sigma_q=\{0,1,\ldots, q-1\}$. For any $n\geq1$, let $\Sigma_q^n=\underbrace{\Sigma_q\times\cdots\times\Sigma_q}_{n}$ be the $n$ copies of $\Sigma_q$, and  $\Sigma_q^\ast=\bigcup_{n=1}^{\infty}\Sigma_q^n$ be the set of all finite words.
        
        \begin{defi}    \label{def:regular_mapping}
        We say that $\tau: \Sigma_q^\ast\to \{-1,0,\ldots,p-2\}$ is a \textit{regular mapping} if
        \begin{enumerate}[(i)]
                \item $\tau(0^n)=0$ for all $n\ge 1;$
                \item $\tau(I)\in (i_n+q \Z)$ for $I=i_1\cdots i_n\in \Sigma^n_q$ for $n\ge 1$;
                \item for any word $I\in \Sigma^\ast_q$, $\tau(I0^l)=0$ for all sufficiently large $l$.
        \end{enumerate}
        \end{defi}
 \noindent Let $\tau: \Sigma_q^\ast\to \{-1,0,\ldots,b-2\}$ be a regular mapping. For any $n\in\N$, there exists a unique $I=i_1i_2\cdots i_N\in\Sigma_q^\ast$ with $i_N\neq0$ such that
\begin{equation}\label{eqadic}
n=i_1+i_2q+\cdots+i_Nq^{N-1}.
\end{equation}
Associated with $\tau$, there is a sequence of integers by $\lambda_0=0$ and
\[\lambda_n=\tau(I|_1)+\tau(I|_2)p+\cdots+\tau(I|_N)p^{N-1}+\sum_{k=N}^\infty\tau(I0^{k-N+1})p^k,\]
from which we note that $\lambda_n$ is uniquely determined by $\tau(I|_1),\tau(I|_2),\ldots,\tau(I|_N)=\tau(I)$. Now by writing $\Lambda(\tau)=\{\lambda_n\}_{n=0}^{\infty}$ and $\ell_n=\#\{k:\tau(I0^k)\neq0,k\geq1\}$ for $n$, we have the following theorem.
\begin{theo}[\cite{DHL13}]\label{D13}
        Let $\tau$ be a regular mapping and let $\Lambda(\tau)$ be defined as above. Then\\
        (i) $\Lambda(\tau)$ is a maximal orthogonal set of $\mu_{p,q}$.\\
        (ii) If $\max_{n\geq1}\{\ell_n\}<\infty$, then $\Lambda(\tau)$ is a spectrum of $\mu_{p,q}$.
\end{theo}
\noindent For simplicity, if $\max_{n\geq1}\{\ell_n\}<\infty$, we call $\Lambda(\tau)$ a {\it regular spectrum} of $\mu_{p,q}$.

By Theorem 3.5 in \cite{DHSW11}, for any orthogonal set of $\mu_{p,q}$, $\dim_{Be}(\Lambda)\leq\frac{\log q}{\log p}$. For the sake of conciseness, we use $s$ to denote this upper bound throughout the article. Among the spectra of $\mu_{p,q}$, the simplest spectrum (also called the {\it canonical spectrum}) is the following set
\[\Lambda_{p,q}=\left\{\sum_{i=1}^{k}a_i p^{i-1}:a_i\in\{0,1,\ldots,q-1\},k\geq1\right\},\]
which attains the maximal Beurling dimension $s$ \cite{DHSW11}.

Compared with the classic case, a natural question is the following:\\
\textbf{Question:} For a maximal orthogonal set of $\mu_{p,q}$, whether we can impose some condition on $D_s^+(\Lambda)$ such that $\Lambda$ is a spectrum?

 In this article, we first investigate the $s$-Beurling density of regular maximal orthogonal sets of $\mu_{p,q}$ and obtain their optimal upper bounds.

\begin{theo}\label{rs}
        For any regular maximal orthogonal sets  $\Lambda$ of $\mu_{p,q}$, $D_s^+(\Lambda)\leq \frac{1}{\left(\frac{q-1}{2(p-1)}\right)^{s}}$. Moreover, the bound is attained by the spectrum $\Lambda_{p,q}$.
\end{theo}      
\begin{rem}
We conjecture that for any maximal orthogonal set $\Lambda$ of $\mu_{p,q}$, we have that $D_s^+(\Lambda)\leq \frac{1}{\left(\frac{q-1}{2(p-1)}\right)^{s}}$ (if it is true, then it is clear that for  any orthogonal set $\Lambda$ of $\mu_{p,q}$, we have that $D_s^+(\Lambda)\leq \frac{1}{\left(\frac{q-1}{2(p-1)}\right)^{s}}$). By Theorem 3.5 in \cite{DHSW11}, we have that for any orthogonal set $\Lambda$ of $\mu_{p,q}$, $D_s^+(\Lambda)<+\infty$, yet it provides a qualitative description rather than quantitative one.
\end{rem}



\begin{rem}
In fact, adopting the same argument as Theorem \ref{rs}, we can show that the following subclass of regular spectra also achieve the largest possible $s$-Beurling density.

\textbf{A}. The authors in
\cite{WW20} gave a new class of spectra for $\mu_{p,q}$ by  using infinite word in $\{-1,1\}^\infty$
acting on $\Lambda_{p,q}$. They proved that for any $\omega=w_1w_2\cdots\in\{-1,1\}^\infty$, the set
\[\Lambda_{w}=\left\{\sum_{i=1}^{k}a_iw_i p^{i-1}:a_i\in\{0,1,\ldots,q-1\},k\geq1\right\}\]
is a spectrum of $\mu_{p,q}$. It is known that for any $\omega=w_1w_2\cdots\in\{-1,1\}^\infty$, we have $\dim_{Be}(\Lambda_w)=s$ \cite{LW}. Based on this result and by using similar way in Theorem \ref{rs}, we can obtain for any $\omega=w_1w_2\cdots\in\{-1,1\}^\infty$, $D_s^+(\Lambda_w)=D_s^+(\Lambda_{p,q})=\frac{1}{\left(\frac{q-1}{2(p-1)}\right)^{s}}.$

\textbf{B}. If the condition $\max_{n\geq1}\{\ell_n\}=0$ (in Theorem \ref{D13}), then $\Lambda(\tau)$ becomes the following set:
\[\Lambda'_{p,q}=\left\{\sum_{i=1}^ka_ip^{i-1}:a_i\in\{0,1,\cdots,q-1\}\pmod{q}\subseteq\{-1,0,1,\cdots,p-2\},k\geq1\right\}.\]
It is easy to see that $\dim_{Be}(\Lambda'_{p,q})=s$. By using similar way in Theorem \ref{rs} again, we have $D_s^+(\Lambda'_{p,q})=\frac{1}{\left(\frac{q-1}{2(p-1)}\right)^{s}}$. \\
On the other hand, as mentioned in the above, there exists a spectrum $\Lambda$ of $\mu_{p,q}$, whose Beurling dimension is zero. This implies that its $s$-Beurling density is zero.
\end{rem}

Finally, we obtain the following.

\begin{theo}\label{thp}
        There exist uncountably many regular spectra $\Lambda$ of $\mu_{p,q}$ such that $D_s^+(\Lambda)$ is positive.
\end{theo}      
        
        
        This paper is organized as follows. In Section 2, we prove Theorems \ref{rs} and \ref{thp}. In Section 3, we conclude with some open questions.

\section{Proofs}
In this section, we prove Theorems \ref{rs} and  \ref{thp}. We start with the following proposition which will be used in Theorem \ref{rs}.

\begin{prop}\label{seq}
        Let $\Lambda$ be a countable subset of $\R$. For any $r>0$,
        \[D_r^+(\Lambda)=2^r\cdot\limsup_{n\to\infty}\sup_{m\in\Z}\frac{\#(\Lambda\cap [m,m+n))}{n^r}.\]
\end{prop}
\proof
For any ball $B(x,h)$, consider the largest interval $[m_1,m_1+n_1)$ contained in $B(x,h)$ (respectively, smallest interval $[m_2,m_2+n_2)$ containing $B(x,h)$) such that $m_1,n_1$ (respectively, $m_2,n_2$) are integers. It is clear that
\[
\frac{\# (\Lambda \cap [m_1,m_1+n_1))}{n_1^r} \cdot \frac{n_1^r}{h^r} \le \frac{\# (\Lambda \cap B(x,r))}{n_1^r} \le \frac{\# (\Lambda \cap [m_2,m_2+n_2))}{n_2^r} \cdot \frac{n_2^r}{h^r}
\]
Now that $\frac{n_1^r}{h^r}$, $\frac{n_2^r}{h^r}$ tends to $2^r$ as $h$ tends to infinity, this implies that
        \[\begin{split}D_r^+(\Lambda)&=\limsup_{h\to\infty}\sup_{x\in\R}\frac{\#(\Lambda\cap B(x,r))}{n^r}=2^r\cdot\limsup_{n\to\infty}\sup_{m\in\Z}\frac{\#(\Lambda\cap [m,m+n))}{n^r}.
        \end{split}\]

\eproof 

First, we prove Theorem \ref{rs}. For any $I=i_1i_2\cdots i_n\in\Sigma_q^n$, denote $|I|=n$ and
\begin{eqnarray*}
    I_{1,k}=
    \begin{cases}
    i_1\cdots i_k, & \hbox{if $k\leq n$;}\\
    i_1\cdots i_n0^{k-n}, & \hbox{if $k>n$.}
    \end{cases}
\end{eqnarray*}
For convenience, we also define $\tau^\ast(I)=\sum_{k=1}^{+\infty}\tau(I_{1,k})p^{k-1}$, so that $\Lambda(\tau)=\{\tau^\ast(I):I\in\Sigma_q^\ast\}$. In the following, we show that the map $\tau^\ast$ is actually a bijection when restricted on the set
\[
\Gamma(\tau)=\{I=i_1 i_2 \cdots i_{|I|}\in\Sigma_q^\ast:i_{|I|} \ne 0\} \cup \{0\},
\]
where $\Gamma(\tau)$ is essentially the set of base-$q$ representations of all nonnegative integers.

\begin{lem}\label{eqv}
The map $\tau^\ast: \Gamma(\tau) \to \Lambda(\tau)$ is a bijection. Furthermore, if $I, J \in \Gamma(\tau)$ are distinct but $I_{1,k}=J_{1,k}$ for some $k \ge 0$, then $|\tau^\ast(I)-\tau^\ast(J)| \ge p^k$.
\end{lem}
\proof
For sujectivity, we show that $\tau^\ast(\Gamma(\tau)) = \Lambda(\tau)$. By definition, $\tau^\ast(\Gamma(\tau)) \subset \Lambda(\tau)$. As for the other inclusion, for any $I = i_1 \cdots i_{|I|} \in \Sigma_q^\ast$ there exists a minimal integer $1 \le N \le |I|$ such that $I=I_{1,N} 0^{|I|-N}$. Then, $I_{1:N} \in \Gamma(\tau)$ and $\tau^\ast(I)=\tau^\ast(I_{1:N})$, which shows $\Lambda(\tau) \subset \tau^\ast(\Gamma(\tau))$. As for injectivity, it suffices to show the proposed inequality. Suppose $I,J\in\Gamma(\tau)$ are distinct and $I_{1,k}=J_{1,k}$. By definition there exists a minimal $M > k$ such that $I_{M,\infty}=J_{M,\infty}$. We then have that $M > k$ and that
\begin{align*}
    & |\tau^\ast(I)-\tau^\ast(J)|=\left|\left(\tau(I_{1,M})-\tau(J_{1,M})\right) p^{M-1}+\sum_{n=k+1}^{M-1} \left(\tau(I_{1,n})-\tau(J_{1,n})\right) p^{n-1}\right| \\
    \ge & \left|\left(\tau(I_{1,M})-\tau(J_{1,M})\right) p^{M-1}\right| - \sum_{n=k+1}^{M-1} \left|\left(\tau(I_{1,n})-\tau(J_{1,n})\right) p^{n-1}\right| \\
    \ge & p^{M-1} - \sum_{n=k+1}^{M-1} (p-1) p^{n-1} = p^k.
\end{align*}
This completes the proof.
\eproof

Now we present the key idea behind Theorem \ref{rs}. Intrinsically, the following lemma illustrates the structure of $\Lambda(\tau)$ and gives a good estimate of the density.
\begin{lem} \label{lem:density_upper_bound}
For any $m\in\Z$ and $k\in\N$,
\[\#(\Lambda(\tau)\cap[m,m+p^k))\leq q^k.\]
\end{lem}
\proof
It is clear by observing that for every $J \in \Sigma_q^k$, the set
\[
S_J:=\{I \in \Gamma(\tau): I_{1:k}=J, \tau^\ast(I) \in [m,m+p^k)\}
\]
contains no more than one element by Lemma \ref{eqv}.
Hence,
\[
\#(\Lambda(\tau)\cap[m,m+p^k))=\sum_{J \in \Sigma_q^k} \# S_J \le q^k.
\]
\eproof

{\noindent{\bf Proof of Theorem \ref{rs}.}}
To prove the first statement of the theorem, let us consider the following density restricted to the interval $[m,m+n)$
\[
\sup_{m\in\Z}\frac{\#(\Lambda\cap [m,m+n))}{n^{s}}.
\]
Representing $n$ in base-$p$ as $n=a_1 + a_2 p + \cdots +a_k p^{k-1}$ ($a_k > 0$), we divide our discussion into the following two cases: \textbf{1.}~$a_i \le q-1$ for $1 \le i \le k$, and \textbf{2.}~there exists $i$ such that $a_i >q-1$. In the former case, one can exploit the fact that $\Lambda$ is a subset of integers and apply Lemma \ref{lem:density_upper_bound} to deduce that
\begin{equation} \label{eq:upper_bound_1}
    \begin{aligned}
        & \sup_{m\in\Z}\frac{\#(\Lambda\cap [m,m+n))}{n^{s}} \le \frac{a_1 + a_2 q + \cdots + a_k q^{k-1}}{(a_1 + a_2 p + \cdots + a_k p^{k-1})^s} \\
        \le & \frac{(q-1) + (q-1) q + \cdots + (q-1) q^{k-1}}{((q-1) + (q-1) p + \cdots + (q-1) p^{k-1})^s},
    \end{aligned}
\end{equation}
where the second inequality follows from monotonicity of the partial derivatives with respect to $a_i$ on the interval $[0,q-1]$. In the latter case, if $i_0$ is the largest index such that $a_{i_0} \ge q$, one can again apply Lemma \ref{lem:density_upper_bound} to show that
\begin{equation} \label{eq:upper_bound_2}
    \begin{aligned}
        & \sup_{m\in\Z}\frac{\#(\Lambda\cap [m,m+n))}{n^{s}} \le \frac{q^{i_0} - 1 + a_{i_0+1} q^{i_0} + \cdots + a_k q^{k-1}}{n^{s}} + \frac{1}{n^{s}} \\
        \le & \frac{(q-1) + (q-1) q + \cdots + (q-1) q^{i_0-1} + a_{i_0+1} q^{i_0} + \cdots + a_k q^{k-1}}{((q-1) + (q-1) p + \cdots + (q-1) p^{i_0-1} + a_{i_0+1} p^{i_0} + \cdots + a_k p^{k-1})^s} + \frac{1}{n^{s}} \\
        \le & \frac{(q-1) + (q-1) q + \cdots + (q-1) q^{k-1}}{((q-1) + (q-1) p + \cdots + (q-1) p^{k-1})^s} + \frac{1}{n^{s}}.
    \end{aligned}
\end{equation}
As a consequence of \eqref{eq:upper_bound_1} and \eqref{eq:upper_bound_2},  $D_s^+(\Lambda)
\leq\frac{1}{(\frac{q-1}{2(p-1)})^s}$ follows immediately by Proposition \ref{seq}.

 Next, we prove the second assertion.
For any $k\geq1$, denote        \[\Lambda_{p,q}^k=\left\{\sum_{i=1}^{k}a_i p^{i-1}:a_j\in\{0,1,\ldots,q-1\}\right\}.\]
Choose $m=0$. Let $n_k=\frac{q-1}{p-1}(p^k-1)$ for $k\geq1$. By Proposition \ref{seq}, we have
\begin{equation} \label{eq:lower_bound}
    \begin{split}
        D_s^+(\Lambda_{p,q})&=2^s\cdot\limsup_{n\to\infty}\sup_{m\in\Z}\frac{\#(\Lambda\cap [m,m+n])}{n^s}\\
        &\geq2^s\limsup_{k\to\infty}\frac{\#(\Lambda\cap [0,n_k])}{n_k^s}\\
        &=2^s\limsup_{k\to\infty}\frac{\#\Lambda_{p,q}^k}{n_k^s}=2^s\limsup_{k\to\infty}\frac{q^k}{(\frac{q-1}{p-1}(p^k-1))^s}\\
        &=\frac{1}{\left(\frac{q-1}{2(p-1)}\right)^{s}}.
    \end{split}
\end{equation}
Combining the first statement of this theorem, we obtain the desired result.

\eproof

Finally, we prove Theorem \ref{thp}.

{\noindent{\bf Proof of Theorem \ref{thp}.}}
We define a mapping by $\tau(0^k)=0$ for $k\ge 1$, and for $n$ and $N$ defined as in \eqref{eqadic}, $\tau(I)=i_N$, and if $n$ is odd, for all $l$,
\begin{equation}\label{mm}
        \tau(I0^l)=0;
\end{equation}
if $n$ is even, let $m_n$ be a strictly increasing sequence of positive integers and  define
\begin{equation}\label{mm}
        \tau(I0^l)=
        \begin{cases}
                0,& \text{if $l\not=m_n;$}\\
                q_{N+m_n},& \text{if $l=m_n.$}
        \end{cases}
\end{equation}
Define $\lambda_0=0$, write $N_n:=N$. When $n\geq1$, if $n$ is odd,
\[\lambda_n=\sum_{j=1}^{N_n}i_jp^{j-1};\]
if $n$ is even,
        \[\lambda_n=\sum_{j=1}^{N_n}i_jp^{j-1}+qp^{N_n+m_n-1}.\]
Then $\Lambda(\tau)=\{0\}\cup\{\lambda_n\}_{n\geq1,n\in2\Z}\cup\{\lambda_n\}_{n\geq1,n\in2\Z+1}=:\Lambda_1\cup\Lambda_2\cup\Lambda_3$.

Consider $\Lambda_2$. Note that for even $n\geq1$, we have $\lambda_n\geq q\cdot p^{m_n+N_n-1}$ and
        \[\begin{split}
                \lambda_n&\leq q\cdot \frac{p^{N_n}-1}{p-1}+q\cdot p^{m_n+N_n-1}\\
                &\leq (q+1)\cdot p^{m_n+N_n-1}.
        \end{split}\]
This implies that
\[\frac{\lambda_{n+2}}{\lambda_n}\geq\frac{q\cdot p^{m_{n+2}+N_{n+2}-1}}{(q+1)\cdot p^{m_n+N_n-1}}\geq \frac{q}{q+1}\cdot p=:b>1.\]
By Proposition 2.2 in  \cite{AL},       $\#(\Lambda_2\cap(x-h,x+h))\leq \log_b (2h+1)$. Consequently, $D_s^+(\Lambda_2)=0$.

Using a similar method as in Theorem \ref{rs}, we have $D_s^+(\Lambda_3)>0$. One can see that  \[D_s^+(\Lambda(\tau))=\max\{D_s^+(\Lambda_1),D_s^+(\Lambda_2),D_s^+(\Lambda_3)\}.\]
Then $D_s^+(\Lambda(\tau))>0$. When $n$ is even larger than $1$. Let $m_n=n^2$ or $n^2+1$. Then we can obtain uncoutably many regular mapping by choosing $m_n$ radomly from the above two choices. Hence, we complete the proof.

\eproof

\section{Further questions}

We end this paper with some further questions. In this paper, we have found a large class of spectrum of $\mu_{p,q}$, which have the biggest $s$-Beurling density. A natural question is the following:

\textbf{Question 3.1.} How to characterize the spectra with maximal s-Beurling density completely?

  Moreover, there is a spectrum of $\mu_{p,q}$, whose $s$-Beurling density is zero. This inspires us to ask the following intermediate  value question with respect to Beurling density.

\textbf{Question 3.2.} Can we always construct a spectrum with any prescibed s-Beurling density from zero to $\frac{1}{\left(\frac{q-1}{2(p-1)}\right)^{s}}$?

On the other hand, we only consider $s$-Beurling density of spectra $\Lambda$ of $\mu_{p,q}$ in this note. We can also consider $r$-Beurling density of spectra $\Lambda$ of $\mu_{p,q}$, where $r\leq s$.

\textbf{Question 3.3.} If $\dim_{Be}(\Lambda)=r~(0<r\leq s)$, find the
upper and lower bounds for the following set
\[\{D_r^+(\Lambda):\dim_{Be}(\Lambda)=r\}.\]

        \subsection*{Acknowledgements}

The authors would like to thank Professor Xinggang He for drawing our attention to these problems and reading the manuscript. The authors would like to thank Professor Meng Wu for his guidance, comments and encouragement when preparing this paper. This project was supported by the Academy of Finland (318217).

\end{document}